\documentclass[11pt, a4paper]{article}
\usepackage[T1]{fontenc}
\usepackage[utf8]{inputenc}
\usepackage{authblk}
\usepackage{amsmath,amsfonts,amssymb,amscd,amsthm,xspace}
\usepackage{graphicx}
\usepackage{vmargin}
\setmarginsrb  { 0.8in}  
               { 0.6in}  
               { 0.8in}  
               { 0.8in}  
               {  20pt}  
               {0.25in}  
               {   9pt}  
               { 0.3in}  

\title{Effective Methods for Solving Band SLEs after Parabolic Nonlinear PDEs}
\author[1]{Milena Veneva \thanks{milena.p.veneva@gmail.com}}
\author[1]{Alexander Ayriyan \thanks{ayriyan@jinr.ru}}
\affil[1]{Joint Institute for Nuclear Research, Laboratory of Information Technologies, Joliot-Curie 6, 141980 Dubna, Moscow region, Russia}

\date{}
\begin{document}
\maketitle
\abstract{%
A class of models of heat transfer processes in a multilayer domain is
considered. The governing equation is a nonlinear heat-transfer equation with
different temperature-dependent densities and thermal coefficients in each
layer. Homogeneous Neumann boundary conditions and ideal contact ones are
applied. A finite difference scheme on a special uneven mesh with a second-order
approximation in the case of a piecewise constant spatial step is built. This
discretization leads to a pentadiagonal system of linear equations~(SLEs) with
a matrix which is neither diagonally dominant, nor positive definite. Two
different methods for solving such a SLE are developed – diagonal
dominantization and symbolic algorithms.
}
\maketitle
%

\section{Introduction and Mathematical Model}
\label{intro}
\begin{figure*}[h]
    \centering
    \includegraphics[width=0.4\textwidth]{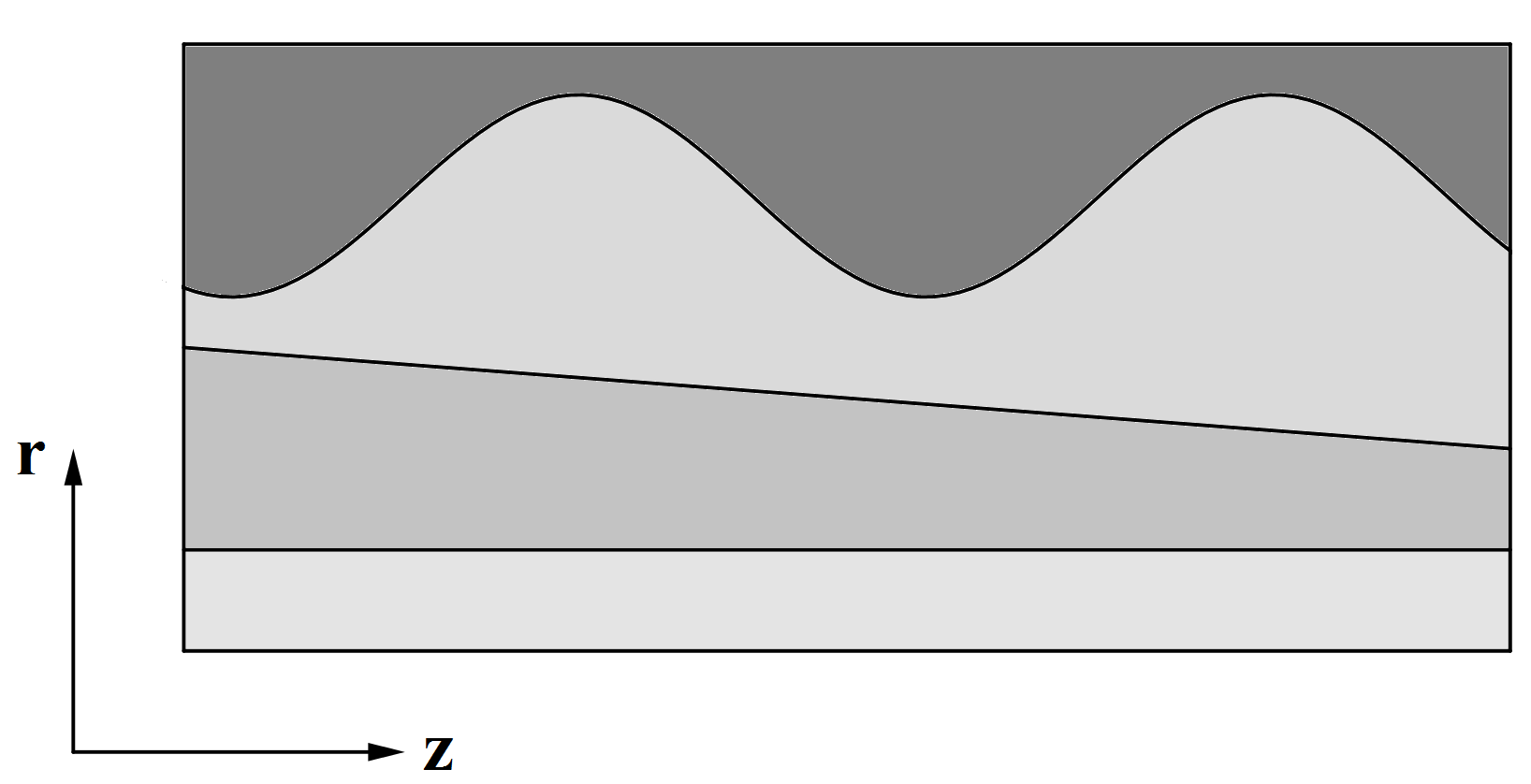}
    \caption{Example of a rectangular domain in cylindrical coordinates~(a
    longitudinal section of a multilayer cylinder), where the thermal
    coefficients are different in all the subdomains and they have a
    discontinuity of the first kind at the borders in-between the subdomains}
    \label{fig_object}
\end{figure*}

In this note, we focus on solving pentadiagonal~(PD) and tridiagonal~(TD)
systems of linear algebraic equations~(SLEs) which are obtained after the
discretization of parabolic nonlinear partial differential equations~(PDEs), using
finite difference methods~(FDM) of second-order approximation. Such a problem was
solved in~\cite{ref_name_1}. There, a finite difference scheme of first-order
approximation was built that leads to a TD SLE with a diagonally dominant
coefficient matrix. The system was solved using the Thomas method (\cite{ref_name_2}).
The following nonlinear model of a cylindrical multilayer structure is considered:
\label{sec-1}
\begin{align}
\label{eq:2_1}
\rho^m (u)\,c_{v}^m (u)\,\frac{\partial u}{\partial t}=\frac{1}{r}\,
\frac{\partial}{\partial r} \left( r\,\lambda^m (u)\,\frac{\partial
u}{\partial r}\right ) + \sum_{\alpha=1}^{n}\frac{\partial}{\partial z_{\alpha}}
\left(\lambda^m (u)\,\frac{\partial u}{\partial z_{\alpha}}\right) + \varphi^m (u)&; \\
\label{eq:2_2}
\frac{\partial u}{\partial r} =0 \qquad \forall \, r \in \{r_{\min},r_{\max}\}&; \\
\label{eq:2_3}
-\lambda^{m}(u)\left.\frac{\partial u}{\partial r}\right|  _{r=r^m_{i^*} -0} = -\lambda^{m+1}(u)\left.\frac{\partial u}{\partial r}\right|  _{r=r^m_{i^*} +0} \quad\textrm{and}\quad
\left.u\right|_{r=r^m_{i^*} -0} = \left.u\right|_{r=r^m_{i^*} +0}&,
\end{align}
where $(r,\vec{z}) \in \Omega\cup\partial\Omega$, $t \geq 0$;
m -- index of the subdomain. Equation~(\ref{eq:2_1}) represents the conservation
of heat inside a multilayer structure. It is an energy equation with conduction
heat transfer, where the densities, the specific heat capacities, and the thermal
conductivities depend on the temperature. For instance, in the two-dimensional
case, Eq.~(\ref{eq:2_1}) could be defined in a domain similar to the one shown
in Fig.~\ref{fig_object}. Homogeneous Neumann boundary conditions~(\ref{eq:2_2})
are applied on the outer boundaries in the radial direction. The ideal contact
conditions~(\ref{eq:2_3}) model the heat flux on the inner boundaries in the
radial direction, where $r^{m}_{i^*}$ denotes the point of discontinuity. The
numerical algorithms for solving the multidimensional governing equation~(GE),
using FDM~(e.g. ADI algorithms (\cite{ref_name_3})), ask for a repeated SLE
solution.
\section{Discretization of the Problem}
\label{sec-2}
Focusing on the radial term of the GE with the assumption that the other terms
will be moved to the right-hand side~(RHS), we consider the following special
mesh with grid points on the inner boundaries:
$\overline{\omega}_{r} = \lbrace (t,r)\left|\right.
t_k = k\, h_t,\,k \in \mathbb{N}_{0};\, r_{0} = r_{\min},\,
r_{i+1} = r_i+h_{i+1},\, i = 0,\ldots,N-2\rbrace.$
A finite difference scheme with a second-order approximation has the following
form~(three-point stencils are taken for the GE and the outer boundary
conditions~(BC), and five-point stencil for the inner BC):
\begin{align}
\label{eq:3_1}
\hat{\rho}_{i}\,\hat{c}_{v_{i}}\frac{\hat{u}_{i}-u_{i}}{\tau}=
\frac{1}{r_i}\,\frac{1}{\hbar_{i}}\left[r_{i+1/2}\,\hat{\lambda}_{i+1/2}\,
\frac{\hat{u}_{i+1}-\hat{u}_{i}}{h_{i+1}}-r_{i-1/2}\,\hat{\lambda}_{i-1/2}\,
\frac{\hat{u}_{i}-\hat{u}_{i-1}}{h_{i}}\right]+\hat{\varphi}_{i}&; \\
\label{eq:3_2}
-\frac{h_{2}\,(2h_{1}+h_{2})\hat{u}_{0}-(h_{1}+h_{2})^{2}\hat{u}_{1}+h_{1}^{2}\hat{u}_{2}}
{h_{1}\,h_{2}\,(h_{1}+h_{2})}= 0&; \\
\label{eq:3_3}
\hat{\lambda}^m_{i^\ast}\,
\frac{h_{i^\ast\!-1}\,(2h_{i^\ast}+h_{i^\ast\!-1})\hat{u}_{i^\ast}-(h_{i^\ast}+h_{i^\ast\!-1})^{2}\hat{u}_{i^\ast\!-1}+h_{i^\ast}^{2}\hat{u}_{i^\ast\!-2}}
{h_{i^\ast}\,h_{i^\ast\!-1}\,(h_{i^\ast}+h_{i^\ast\!-1})} =&\\\nonumber =-\hat{\lambda}^{m+1}_{i^\ast}\,\frac{h_{i^\ast\!+2}\,(2h_{i^\ast\!+1}+h_{i^\ast\!+2})\hat{u}_{i^\ast}-(h_{i^\ast\!+1}+h_{i^\ast\!+2})^{2}\hat{u}_{i^\ast\!+1}+h_{i^\ast\!+1}^{2}\hat{u}_{i^\ast\!+2}}
{h_{i^\ast\!+1}\,h_{i^\ast\!+2}\,(h_{i^\ast\!+1}+h_{i^\ast\!+2})}&; \\
\label{eq:3_5}
\frac{h_{N-2}\,(2h_{N-1}+h_{N-2})\hat{u}_{N-1}-(h_{N-1}+h_{N-2})^{2}\hat{u}_{N-2}+h_{N-1}^{2}\hat{u}_{N-3}}
{h_{N-1}\,h_{N-2}\,(h_{N-1}+h_{N-2})}= 0&,
\end{align}
where
\begin{equation*}
\hat{\lambda}_{i\pm 1/2} = \lambda\left(\frac{\hat{u}_{i}+\hat{u}_{i\pm 1}}{2}\right), \qquad\hbar_{i}=\frac{h_{i+1}+h_{i}}{2},\qquad r_{i\pm 1/2} = \frac{r_{i}+r_{i\pm 1}}{2}.
\end{equation*}
The matrix form of the considered system is: $A\vec{\hat{u}} = \vec{\varphi}(\hat{u})$,
where $A$ is a PD sparse matrix which does not have any special properties, e.g.
diagonally dominance or positive definiteness. In order to preserve the band
structure of the matrix, we cannot use the Gaussian elimination with pivoting.
Within this, we could obtain division by zero at some point of the procedure, which
is going to make our algorithm unstable. For that reason, we alter the initial PD
matrix by adding the minimum values to the diagonal elements so as to transform the
matrix into a weakly diagonally dominant one:
\begin{align*}
A_{DD}\,\vec{\hat{u}} = \vec{\varphi}(\hat{u}) + P\,\vec{\hat{u}}&, \qquad\textrm{where}\qquad A_{DD} = A + P;\\
P = \textrm{diag}\left(2h_{1}^2, \delta_{i^\ast\!,j}\,p_{i^\ast\!,j}, 2h_{N-1}^2\right)&, \qquad\textrm{where}
\qquad p_{i^\ast\!,i^\ast} = \frac{2\lambda_{i^{\ast}}^{m}\,h_{i^\ast}}{h_{i^\ast\!-1}\,(h_{i^\ast}+h_{i^\ast\!-1})}+\frac{2\lambda_{i^{\ast}}^{m+1}\,h_{i^\ast\!+1}}{h_{i^\ast\!+2}\,(h_{i^\ast\!+1}+h_{i^\ast\!+2})}.
\end{align*}
The Gaussian elimination with pivoting~(the procedure we use to transform
the~PD matrix into a~TD one) does not preserve the diagonal dominance of the
matrix. The use of the Gaussian elimination to the initial matrix
$A$~(not $A_{DD}$) yields a transformed matrix $\tilde{A}$. Then,
the diagonal dominantization method is used in order to transform the obtained
TD matrix $\tilde{A}$ into a diagonally dominant one. To that purpose,
the nondiagonal elements are added to the diagonal ones:
\begin{align*}
\tilde{A}_{DD}\,\vec{\hat{u}} = \vec{\varphi}(\hat{u}) + \tilde{P}\,\vec{\hat{u}}&, \qquad\textrm{where}\qquad\tilde{A}_{DD} = \tilde{A} + \tilde{P};\\
\tilde{P} = \textrm{diag}\left(|\tilde{A}_{0,1}|, \delta_{i^\ast\!,j}\,\tilde{p}_{i^\ast\!,j}, |\tilde{A}_{N-2,N-1}|\right)&, \qquad\textrm{where}\qquad\tilde{p}_{i^\ast\!,i^\ast} = \sum\limits_{\scriptscriptstyle{\beta \in \{-1,1\}}}|\tilde{A}_{i^\ast\!,i^\ast\!+\beta}|.
\end{align*}
\section{Numerical and Symbolic Algorithms}
\label{sec-4}
Two different approaches for solving the SLE are considered -- numerical and
symbolic. The complexity of all the suggested numerical algorithms is $O(N)$.
Since it is unknown what stands behind the symbolic library, evaluating the
complexity of the symbolic algorithms is a very complicated task.

\textbf{Numerical Algorithms.} Two different numerical algorithms are applied
to the system with a~PD matrix. Both of them are based on LU~decomposition.
The first one~(\cite{ref_name_4}) is intended for a dense~PD matrix. In the
case of the considered problem, the matrix is sparse. For that reason, a
modified algorithm is built. The main idea is that after the mesh was defined,
the indexes of the discontinuity points are known. Since these indexes
coincide with those of the matrix rows which are not sparse, we can
reference them to the algorithm and conduct the full calculation only for
them. For the rest of the rows, the algorithm is reduced to a problem similar
to solving a system with~TD matrix. This way, the complexity of the algorithm
is decreased but at the cost of additional $N+2$ check-ups for the non-sparse
rows. In the case of the~TD matrix, the system is solved using the Thomas method.

\textbf{Symbolic Algorithms.} The symbolic algorithm in the case of a PD matrix
is also based on LU decomposition~(\cite{ref_name_4}). For the~TD matrix, a
symbolic version of Thomas method is created. As it is known, Thomas method is not
suitable for nondiagonally dominant matrices (\cite{ref_name_2}). In order to
cope with this problem, in the case of a zero quotient of two subsequent leading
principal minors, a symbolic zero is assigned instead and the calculations are
continued. At the end of the algorithm, this symbolic zero is substituted with zero.
The same approach is suggested in~\cite{ref_name_5}.

\section{Implementation and Results}
\label{sec-5}
All the numerical algorithms are implemented using \texttt{C++}. The matrix needs
to be nonsingular and diagonally dominant so as the methods to be stable. Two
symbolic algorithms are implemented, using the \texttt{GiNaC}
library~(version 1.7.2)~(\cite{ref_name_6}) of \texttt{C++}.
The symbolic algorithms require the matrix to be nonsingular only.
In Table~\ref{fig:tab-1} one can find the wall-clock time results from the conducted
experiments. Since the largest supported precision in the
\texttt{GiNaC} library is \textbf{double}, during all the experiments double data
type is used. The notation is as follows: \textbf{NPDM} stands
for numerical~PD method, \textbf{MNPDM} -- modified numerical~PD method,
\textbf{SPDM} -- symbolic~PD method, \textbf{NTDM} -- numerical~TD method,
\textbf{STDM} -- symbolic~TD method. The achieved accuracy is summarized, using
infinity norm. On the penultimate row of the table, one can find the complexity of
all the considered methods. On the last row, the characteristics of the computer
which is used are described. The number of considered discontinuity points is $K = 11$.
\renewcommand{\arraystretch}{1.2}
\begin{table}[!htb]
\caption{Results from solving SLE}
\label{fig:tab-1}
\centerline{\begin{tabular}{cccccc}
\hline\hline
& \multicolumn{5}{c}{Wall-clock time [s]}  \\\hline\hline
$N$ & NPDM & MNPDM & SPDM & NTDM & STDM \\\hline
$10^{3}$ & 0.000036 & 0.000034 & \quad\;\,0.088669 & 0.000021 & \quad\,0.043690 \\
$10^{4}$ & 0.000403 & 0.000373 & \quad\;\,8.467241 & 0.000245 & \quad\,2.971745 \\
$10^{5}$ & 0.004709 & 0.003916 & 3547.020851 & 0.002416 & 799.533587 \\
$10^{8}$ & 3.159357 & 2.682258 & -- & 1.652945 & -- \\
\hline\hline
$\max\limits_{\textrm{\scriptsize{N}}}\|y - \bar{y}\|_{\infty}$ & $2.22\times 10^{-16}$ & $2.22\times 10^{-16}$ & 0 & $2.22\times 10^{-16}$ & 0 \\
\hline\hline
Complexity: & $19N - 29$ & $13N+7K-14$ & -- & $9N+2$ & -- \\\hline\hline
\multicolumn{6}{l}{OS: Fedora 25; Processor: Intel Core i7-6700 (3.40 GHz); Compiler: GCC 6.3.1 (-O0).}\\
\hline\hline
\end{tabular}}
\end{table}

\section{Discussion and Conclusions}
\label{sec-6}
A nonlinear heat transfer equation in a multilayer domain was considered. The
suggested discretization scheme always has a first-order approximation. In the
case of piecewise constant thermal conductivities or when
$\|\hat{\lambda}_{i+1/2}-\hat{\lambda}_{i-1/2}\|_{\infty} = O(\|h_{i}\|_{\infty})$,
the order of approximation is going to be second. Focusing on the radial term, a
SLE with a~PD matrix was obtained. Then, applying Gaussian elimination, a~TD matrix
was derived. For both these matrices, a diagonal dominantization procedure was
suggested. This approach ensures the stability of the suggested methods. A modified
version of the numerical method for solving a SLE with a~PD matrix was built. Since
the complexity of this method is lower than the complexity of the general
algorithm~(usually $K\ll N$), better computational time was achieved. The fastest
numerical algorithm was found to come from the Thomas method. All the experiments
gave an accuracy of an order of magnitude of $10^{-16}$. As a next step symbolic
algorithms were used. They do not require the matrices to be of a special form and
are exact. However, they are not comparable with the numerical algorithms with
respect to the required time in the case of a numerical solving of the heat equation
when one needs to solve the SLE many times. On the other hand, these symbolic methods
are not as restrictive as the numerical ones when it comes to the matrix properties.
Another upside of the symbolic algorithms is that in the case of a piecewise linear
equation, they do not add nonlinearity to the RHS of the system and hence, there is
no need of iterations for the time step to be executed. In future, the approach
suggested in this note will be investigated in detail.

\section*{Acknowledgements}
The authors want to express their gratitude to the Summer Student Program at JINR, and J\'{a}n Bu\v{s}a Jr. (JINR). A.~Ayriyan thanks the JINR grant No.~17-602-01.


\end{document}